\documentclass[12pt]{article}
\usepackage{amsmath,amsfonts,amssymb}
\usepackage{color}
\newtheorem{theorem}{\hspace*{\parindent}Theorem}
\newtheorem{lemma}{\hspace*{\parindent}Lemma}

\newcounter{theremark}

\newcommand{\rr}{\mathbb R^d}

\title{Optimal Green energy points on the circles in $d$-space}
\author{V.N.\:Dubinin$^{\rm a,b}$\footnote{Corresponding author. E-mail addresses: \emph{dubinin@iam.dvo.ru} (V.N.\:Dubinin),  \emph{pril-elena@yandex.ru} (E.G.\:Prilepkina)}~~and E.G.\:Prilepkina$^{\rm b}$
\\[10pt]
\small{\textit{$\phantom{1}^a$Far Eastern Federal University
(FEFU), 8, Sukhanova Street, Vladivostok, 690950,
Russia}}\\\small{\textit{$\phantom{1}^b$Institute of Applied
Mathematics, FEBRAS, 7, Radio Street, Vladivostok, 690041,
Russia}}}
\date{}

\begin{document}
\maketitle

\begin{abstract}
We give two precise estimates for the Green energy of a discrete
charge, concentrated in the points on the circles, with respect to
the concentric rotation domain in the d-dimensional Euclidean
space, $d>2$.The proof is based on the application of a
dissymmetrization, extremal metrics approach and an asymptotic
formula for the condenser capacities in the case when some of its
plates contract to given points.
\end{abstract}

\bigskip


\bigskip


\bigskip

\section{Introduction and statement of results}

The Riesz $s$--energy ($s\neq 0$) of $n$ points $z_1,\ldots,z_n$
of the complex plane is defined by

$$ \sum\limits_{k=1}^n\sum
\limits_{{l=1} \atop {l\not= k}}^n|z_k-z_l|^{-s}.$$ It can be
shown using the classical T\.{o}th's result \cite[p.155]{Fejes}
and a convexity argument that for $s\geq -1$ and each $n\geq 2$,
the $n$--th roots of unity $z_k^*=\exp\{2\pi i(k-1)/n\}$,
$k=1,\ldots,n,$ form minimal $n$--point $s$--energy configuration
for the unit circle $|z|=1,$
\begin{equation}\label{eq:ineq1}
\sum\limits_{k=1}^n\sum \limits_{{l=1} \atop {l\not=
k}}^n|z_k-z_l|^{-s}\geq \sum\limits_{k=1}^n\sum \limits_{{l=1}
\atop {l\not= k}}^n|z_k^*-z_l^*|^{-s}.
\end{equation}
Various sophisticated problems related to the optimality of the
Riesz $s$-energy for different values of $s$ and for the points
$z_k$ lying in the plane sets or in $\mathbb R^d$ have been
treated in a number of papers (see, for instance,
\cite{Brauchart}-\cite{Brauchart2}, and references therein). In
this note we consider the Green energy with respect to a rotation
domain in $\mathbb R^d$,  $d\ge3$,  of a discrete charge
concentrated in the points of certain circles. Limit cases of the
optimal properties of this energy lead to inequalities for the
Riesz $s$-energy for $s=d-2$.  Unlike previous works we study the
charges of the opposite signs located on a collection of circles.
We now turn to the precise formulations. In what follows  $\rr$ is
a  $d$-dimensional Euclidean space with the usual norm
$\|\cdot\|$, with points $\mathbf{x}=(x_1,\dots,\,x_d)$, $d\ge3$.
A domain  $B$ in $\rr$ is admissible if it has the Green function
for the Laplace operator vanishing at the points of the boundary
$\partial B$ of the domain  $B$.  This Green function with pole at
the point  $\mathbf x_0\in B$ will be denoted by  $g_{\small
B}(\mathbf x, \mathbf x_0)$. In the neighborhood of $\mathbf x_0$
the following expansion holds
$$
g_{\small B}(\mathbf x,\mathbf x_0)=\lambda_d(\|\mathbf x-\mathbf
x_0\|^{2-d}-(r(B,\mathbf x_0))^{2-d}+o(1)),\qquad \mathbf x\to
\mathbf x_0,
$$
where $\lambda_d=((d-2)\omega_{d-1})^{-1}$,
$\omega_{d-1}=2\pi^{d/2}/\Gamma(d/2)$ is the surface measure of
the unit hyper-sphere. In all points of $B$ different from the
pole $\mathbf x_0$, the Green function is harmonic, that is
$\triangle g_{\small B}(\mathbf x,\mathbf x_0)=0$. The quantity
$r(B,\mathbf x_0)$ is known as the harmonic radius of the domain
$B$ with respect to the points $\mathbf x_0$ \cite{Ban}.

Denote by  $J$ $(d-2)$--dimensional plane $\{\mathbf
x\in\rr:\mathbf x=(0,0,x_3,\ldots, x_d)\}$. We will need the
cylindrical coordinates $(r,\theta, \mathbf{x}')$ of the point
$\mathbf{x}=(x_1,\dots,\,x_d)$ in $\rr$,  related to the Cartesian
coordinates by $x_1=r \cos \theta,$ $x_2=r \sin \theta,$
$\mathbf{x}'\in J$. A domain $B\subset \rr$ will be called
\emph{the rotation domain} (with respect to the axis $J$), if for
any point $(r,\theta,\mathbf{x}')\in B$ and any $\varphi$ the
point $(r,\varphi,\mathbf{x}')$ belongs to $B$.

Suppose that $B$ is an admissible rotation domain and let
$\Omega=\{S\}$ be the collection comprising a finite number of
distinct circles $S$ of the form $S=\{(r_0,\theta,\mathbf
x'_0):0\leq\theta\leq 2\pi\}$ lying in the domain $B$ (here
$r_0>0$ and $\mathbf x'_0\in J$ are assumed to be fixed). For
arbitrary real numbers $\theta_j,$ $j=0,\ldots,$ \makebox{$m-1,$}
\begin{equation*}
0\leq \theta_0<\theta_1<\ldots<\theta_{m-1}<2\pi,
\end{equation*}
denote by $X=\{\mathbf x_k\}_{k=1}^n$ the collection of all
distinct points of $B$ at which the circles  from  $\Omega$
intersect the half-planes
\begin{equation*}
L_j=\{(r,\theta,\mathbf x'):\theta=\theta_j\}, \ j=0,\ldots,m-1.
\end{equation*}
Let $\Delta=\{\delta_k\}_{k=1}^n$ be an arbitrary discrete charge
(a collection of real numbers), having the value $\delta_k$ at the
point $\mathbf x_k$, $k=1,\ldots,n$. \emph{The Green energy} of
this charge with respect to the domain $B$ is defined by
\cite{Lan}
$$
E(X,\Delta,B)=\sum\limits_{k=1}^n\sum \limits_{{l=1} \atop {l\not=
k}}^n \delta_k\delta_l g_{\small B}(\mathbf x_k,\mathbf x_l).
$$
Define also  $X^*=\{\mathbf x^*_k\}_{k=1}^n$ -- the collection of
points at which the circles from  $\Omega$ intersect the
half-planes
\begin{equation*}
\ L_j^*=\{(r,\theta,\mathbf x'): \theta=2\pi j/m\}, \
j=0,\ldots,m-1.
\end{equation*}

\begin{theorem}
Suppose that the charge $\Delta=\{\delta_k\}_{k=1}^n$ takes equal
values $\delta_k=\delta_l$ at the points $\mathbf x_k$ and
$\mathbf x_l$ from the collection $X$ that lie on the same circle
from $\Omega$ and, furthermore, that the points $\mathbf x_k\in X$
and $\mathbf x^*_k\in X^*$ lie on the same circle from $\Omega$,
$k=1,\ldots,n$. Then
$$
E(X,\Delta,B)\geq E(X^*,\Delta,B).
$$
\end{theorem}

The following proposition asserts that under certain conditions
the symmetric configuration has the maximal energy.

\begin{theorem}
Suppose that the domain $B$ and the collections $\Omega$, $X$ and
$X^*$ are as defined above while $m$ is an even number. Assume
further that the charge $\Delta=\{\delta_k\}_{k=1}^n$ takes the
values of equal moduli $|\delta_k|=|\delta_l|$ at the points
$\mathbf x_k$ and $\mathbf x_l$ belonging to $X$ and lying on the
same circle from $\Omega$ and, moreover $\delta_k<0$ if the point
$\mathbf x_k$ belongs to one of the a half-planes $L_{2 p-1}$,
$1\leq p\leq m/2$, otherwise $\delta_k>0$, $k=1,\ldots,n$. Then
$$
E(X,\Delta,B)\leq E(X^*,\Delta,B),
$$
where the points of the collection $X^*$ are numbered as follows:
if $\mathbf x_k^*\in X^*$ lies at the intersection of a circle $S$
from  $\Omega$ and a half-plane $L_j^*$, then the corresponding
point $\mathbf x_k\in X$ must lie at the intersection of the
circle $S$ and the half-plane $L_j$, $k=1,\ldots,n$,  $0\leq j\leq
m-1$.
\end{theorem}

Note that using the symmetry principle for harmonic functions it
is not difficult to establish \cite{Ban} that the Green function
of the ball $B(0,t)=\{\mathbf x \in \rr:\|\mathbf x\|<t\}$ with
the pole at the point $\mathbf x_0\in B(0,t)$ takes the form
$$
g_{B(0,t)}(\mathbf x, \mathbf x_0)= \lambda_d\left(\|\mathbf x-
\mathbf x_0\|^{2-d}-\left\| \frac{\|\mathbf x_0\|}{t}\mathbf
x-\frac{t }{\|\mathbf x_0\|}\mathbf x_0 \right\|^{2-d}\right).
$$

Placing the points from  $X$ into a sufficiently large ball
$B(0,t)$ and letting $t\to\infty$, from Theorem~1,2 we deduce
inequalities for the  Riesz $(d-2)$--energy.  In particular,
Theorem~1 leads to inequality \eqref{eq:ineq1} for $s=d-2$, while
Theorem~2 yields
\begin{equation*}\sum\limits_{k=1}^{2n}\sum
\limits_{{l=1} \atop {l\not= k}}^{2n}\frac{(-1)^{k+l}}{| z_k-
z_l|^{d-2}}\leq \sum\limits_{k=1}^{2n}\sum \limits_{{l=1} \atop
{l\not= k}}^{2n}\frac{(-1)^{k+l}}{| z_k^*-z_l^*|^{d-2}},
\end{equation*}
where $z_k$, $k=1,\ldots,2n,$ are located on the circle $|z|=1$ in
the ascending order of the index $k$ and $z_k^*=\exp\{\pi
i(k-1)/n\}$,  $k=1,\ldots,2n$.

The proofs of Theorems~1,2 hinge on the theory of condenser
capacity and dissymmetrization  \cite{Dub3}, \cite{DubininBook}.
These proofs are related conceptually with the solutions of the
so-called extremal decomposition problems \cite{Dubinin1},
\cite{DubininPril}, \cite{GKP}, \cite{KalmPrilEn}. In the recent
paper \cite{Dubinin1}, analogues of Theorems~1,2 for the case of
the plane and one circle and a concentric ring have been
presented. The proof of Theorem~1 of this paper follows the same
line of argument as the one presented in \cite{Dubinin1} with
modifications related to the use of dissymmetrization \cite{Dub3}
and the asymptotic formula for the capacity of the spacial rather
than plane condenser \cite{DubininPril}.  The proof of an analogue
of Theorem~2 for the plane case \cite{Dubinin1}  is based on the
radial averaging transformation and conformal mapping.  This
method is not applicable in the Euclidean space due to absence of
the suitable conformal mappings.  Therefore in order to
demonstrate Theorem~2 we resort to the moduli of the families of
curves (see, for  instance, \cite{Fuglede}, \cite{Alf},
\cite{Otsuka}). The idea behind this approach goes back to the
proof of Theorem~4 from \cite{DubininPril}. Our results, as well
as their proofs, can be carried over to the discrete energy with
the Robin function kernel \cite{Duren} (of the domain $B$ with
respect to a part of the boundary) in place of the Green function
kernel. The next section is of an auxiliary nature.

\section{Preliminaries}

Suppose $B$ is an admissible domain in the space $\rr$, $d>2$;
$X=\{\mathbf x_k\}_{k=1}^n$ is a collection of distinct points in
$B$; $\Lambda =\{\sigma_k\}_{k=1}^n$  is a collection of
non-vanishing real numbers; $\Psi=\{\mu_k\}_{k=1}^n$ is a
collection of positive numbers $\mu_k$. Denote by $E(\mathbf
a,t)=\{\mathbf x \in \rr:\|\mathbf x-\mathbf a\|\leq t\}$ the
closed ball of radius $t$ centered at $\mathbf a$. For
sufficiently small $t>0$ introduce  ''the generalized'' condenser
as the ordered collection
$$
C(t;B,X,\Lambda,\Psi) =\{\bar\rr\setminus B, E(\mathbf x_1, \mu_1
t),...,E(\mathbf x_n, \mu_n t)\}
$$
with pre-assigned values $0, \sigma_1,...,\sigma_n$, respectively
\cite{DubininPril}. Similarly to the usual condensers,  defined
the capacity (or 2-capacity) of the condenser
$C(t;B,X,\Lambda,\Psi)$ by
$$
{\mathrm{cap}}\,C(t;B,X,\Lambda,\Psi)=\inf\int_{\rr}|\nabla
v|^2dx,
$$
where the infimum is taken over all functions $v:\bar\rr\to
\mathbb{R}$ from $C^{\infty}(\rr)$, vanishing in a neighborhood of
the set ${\bar \rr}\setminus B$ and equalling to $\sigma_l$ in a
neighborhood $E(\mathbf x_l,\mu_l r),\ l=1,...,n$. The condenser
modulus $|C(t;B,X,\Lambda,\Psi)|$ is reciprocal to the capacity of
$C(t;B,X,\Lambda,\Psi)$:
$$
|C(t;B,X,\Lambda,\Psi)|=\left({\mathrm{cap}}\,C(t;B,X,\Lambda,\Psi)\right)^{-1}.
$$

\begin{lemma}\label{asimptotika}\emph{\cite[Theorem 1]{DubininPril}.}
The following asymptotic formula holds as $t\to 0$:
\begin{equation}\label{eq:asimptotika}
|C(t;B,X, \Lambda,\Psi)|=\nu\lambda_{d}t^{2-d}-\lambda_d
\nu^2\sum\limits_{k=1}^n \nu_k^2r(B,\mathbf x_k)^{2-d}+
\nu^2\sum\limits_{k=1}^{n}\sum\limits_{{l=1} \atop {l\not=
k}}^{n}\nu_l\nu_k g_{\small B}(\mathbf x_l,\mathbf x_k)+o(1),
\end{equation}
where $\nu_k=\sigma_k\mu_k^{d-2},$ $\nu=\left(\sum\limits_{k=1}^n
\sigma^{2}_{k} \mu^{d-2}_{k}\right)^{-1},$ $k=1,...,n.$
\end{lemma}

Let $\Gamma$ be a family of curves in $\rr$. We will assume that
each curve $\gamma\in\Gamma$ is a union of a countable number of
open arcs, closed arcs or closed curves and is locally
rectifiable. {\it $2$-modulus} or just  {\it modulus} of the
family $\Gamma$ is defined as the quantity
$$
M(\Gamma)=\inf\int_{\rr}\rho^2 dx,
$$
where infumum is taken over all Borel functions
$\rho:\rr\rightarrow[0,\infty]$ such that
$\int_{\gamma}{\rho}ds\ge1$ holds for each curve $\gamma\in\Gamma$
\cite{Otsuka}.  It is said that the family $\Gamma_2$ is minorized
by the family $\Gamma_1$, if each curve $\gamma\in\Gamma_2$ has a
sub-curve belonging to $\Gamma_1$.  The families $\Gamma_1,
\Gamma_2,\ldots$ are called separated if there exist disjoint
Borel sets $E_i$ in $\rr$, such that $\int_{\gamma}\chi_i ds=0$
for any curve $\gamma\in\Gamma_i$, where $\chi_i$ is the
characteristic function of $\rr\setminus E_i$. If $\Gamma_1,
\Gamma_2,\ldots$ are separated families and $\Gamma_i$ is
minorized by $\Gamma$, $i=1,2,\dots,$, then
\begin{equation}\label{prop1}
M(\Gamma)\geq \sum_{i=1}^{\infty}M(\Gamma_i).
\end{equation}
If, on the contrary, $\Gamma$ is minorized by $\Gamma_i$,
$i=1,2,\dots,$ and $\Gamma_1, \Gamma_2,\ldots$ are separated
families, then
\begin{equation}\label{prop2}
M(\Gamma)^{-1}\geq \sum_{i=1}^{\infty}M(\Gamma_i)^{-1}.
\end{equation}

It is easy to see that the capacity of the condenser
$C(t;B,X,\Lambda,\Psi)$ under the choice $\sigma_k=1$,
$k=1,\ldots,n$, coincides with the capacity of the condenser with
two plates $E(\mathbf x_1, \mu_1 t)\cup E(\mathbf x_2, \mu_2
t)\ldots\cup E(\mathbf x_n, \mu_n t)$ and $\overline{\mathbb
R}^d\setminus B$ (for the definition of the condenser capacity
see, for instance, in \cite{Dub3}, \cite{Hesse}).  Therefore, the
following lemma holds true.

\begin{lemma}\label{lemma2}\emph{\cite{Hesse}}.
Let $\sigma_1=...=\sigma_n=1$ or $\sigma_1=...=\sigma_n=-1$,
$\Lambda =\{\sigma_k\}_{k=1}^n$, $B,X,\Psi$ as defined above,
$\Gamma(t;B,X,\Psi)$ is the family of continuous curves in $B$
connecting the set $E(\mathbf x_1, \mu_1 t)\cup E(\mathbf x_2,
\mu_2 t)\ldots \cup E(\mathbf x_n, \mu_n t)$ with the boundary
$\partial B$ of the domain $B$. Then
$$
{\mathrm{cap}}\,C(t;B,X,\Lambda,\Psi)=M(\Gamma(t;B,X,\Psi)).
$$
\end{lemma}

We will further need the definition of dissymmetrization in
Euclidean space \cite{Dub3}. Denote by $\Phi$ the group of
reflections in $\overline{\mathbb R}^d$ with respect hyper-planes
of the form $\{(r,\theta,\mathbf x'):\theta=\pi k/m,\ \text{or}\
\theta=\pi+\pi k/m\}, \ k=1,\ldots,m$. Next we introduce a
symmetric structure $\{P_i\}_{i=1}^N$ in $\overline{\mathbb R}^d$
as the  collection of closed angles $P_i=\{(r,\theta,\mathbf x'):
\theta_{i1}\leq \theta\leq\theta_{i2},\, 0\leq r\leq\infty\}$,
$i=1,\dots,N$, satisfying the conditions:
$$
\begin{array}{ll}
\text{aP)} & \bigcup_{i=1}^N P_i = \overline{\mathbb R}^d, \ \sum_{i=1}^N (\theta_{i2}-\theta_{i1})=2\pi,\\
\text{bP)} & \{\phi(P_i)\}_{i=1}^N=\{P_i\}_{i=1}^N \ \text{for any
isometry} \ \phi\in \Phi.
\end{array}
$$

The family of rotations $\{\lambda_i\}_{i=1}^N$ of the form
$\lambda_i(r,\theta,\mathbf x')=(r,\theta+\varphi_i,\mathbf x')$,
$i=1,\ldots,N$, will be called the dissymmetrization of the
symmetric structure $\{P_i\}_{i=1}^N$, if the images
$S_i=\lambda_i(P_i)$ satisfy the following conditions:

\medskip

aS) $\bigcup_{i=1}^{N}S_i=\overline{\mathbb R}^d$,

\medskip

bS) for any non-empty intersection $S_i\cap S_j$, $i,j=1,\dots,N$,
there exists an isometry $\phi~\in~\Phi$, such that
$\phi(\lambda_i^{-1}(S_i\cap S_j))=\lambda_j^{-1}(S_i\cap S_j)$.

\medskip

For an arbitrary set $A$ in $\overline{\mathbb R}^d$ introduce the
notation ${\rm Dis}\,A = \bigcup\limits_{i=1}^{N}\lambda_i(A\cap
P_i)$. A characteristic feature of a rotation domain $B$ is the
fact that such domain is invariant with respect to any
dissymmetrizaton ${\rm Dis}\,B=B$.

\begin{lemma}\label{LemmaDis}\emph{\cite{Dub3}.}
Let the numbers $\theta_j$, $j=0,\ldots,m$, satisfy
$0\le\theta_0<\theta_1<\ldots<\theta_{m-1}<2\pi$,
$\theta_m=\theta_0+2\pi$, and suppose that
$L_j=\{(r,\theta,\mathbf x'):\theta=\theta_j\}$,
$L_j^*=\{(r,\theta,\mathbf x'): \theta=2\pi j/m\},
j=0,\ldots,m-1$. Then there exists a symmetric structure
$\{P_i\}_{i=1}^N$, $N\geq m$, and a dissymmetrization
$\{\lambda_i\}_{i=1}^N$,  such that ${\rm Dis}\,L_j^*=L_j,$
$j=0,\ldots,m-1$, and each half-plane $L_j^*$ is the bisector of a
dihedral angle $P_i$ of size $\psi$, where
$$
\psi=\min\limits_{i=1,\ldots, m}(\theta_i-\theta_{i-1}).
$$
\end{lemma}

The condenser $C(t;B,X,\Lambda,\Psi)$ will be called symmetric
with respect to the group $\Phi$,  if $B$ is a rotation domain and
for any $k$, $k=1,\ldots,n$, and any isometry $\phi\in\Phi$ we
have $\phi(\mathbf x_k)\in X$ and $\sigma_k=\sigma_l$,
$\mu_k=\mu_l$ in the case  $\phi(\mathbf x_k)=\mathbf x_l$. The
result of dissymmetrization of a symmetric condenser
$C(t;B,X,\Lambda,\Psi)$ is defined to be the condenser ${\rm
Dis}\,C(t;B,X,\Lambda,\Psi)=C(t;B,\{{\rm Dis}\,\mathbf
x_k\}_{k=1}^n,\Lambda,\Psi)$.

\begin{lemma}\label{lemma4}
If the condenser $C(t;B,X,\Lambda,\Psi)$ is symmetric with respect
to the group $\Phi$, then for sufficiently small $t$ the following
inequality holds
$$
|C(t;B,X,\Lambda,\Psi)|\leq| {\rm Dis}\,C(t;B,X,\Lambda,\Psi)|.
$$
\end{lemma}

The proof of this claim is essentially the same as the proof of a
similar statement in  \cite[Theorem~4.14]{DubininBook}. A
particular case has been considered in \cite[Theorem~5]{Dub3}.

\section{Proofs of Theorems}

We will start with the proof of Theorem~1.  Suppose the domain $B$
and the collection $X=\{\mathbf x_k\}_{k=1}^n$, $X^*=\{\mathbf
x_k^*\}_{k=1}^n$, $\Delta=\{\delta_k\}_{k=1}^n$ as in Theorem~1.
We can assume that $\delta_k\neq 0$, $k=1,\ldots,n$. Put
$\sigma_k=\text{sgn}\,\delta_k$, $\mu_k=|\delta_k|^{1/(d-2)}$,
$k=1,\ldots,n$, $\Lambda =\{\sigma_k\}_{k=1}^n$,
$\Psi=\{\mu_k\}_{k=1}^n$.  Note that the condenser
$C(t;B,X^*,\Lambda,\Psi)$ is symmetric with respect to the group
$\Phi$. Apply dissymmetrization from Lemma~\ref{LemmaDis} to the
condenser $C(t;B,X^*,\Lambda,\Psi)$.  The result of
dissymmetrization of this condenser for small $t$ is the condenser
$C(t;B,X,\Lambda,\Psi)$.  According to Lemma~\ref{lemma4}
$$
|C(t;B,X^*,\Lambda,\Psi)|\leq |C(t;B,X,\Lambda,\Psi)|.
$$
Applying the asymptotic formula \eqref{eq:asimptotika}, we obtain
\begin{multline}\label{as2}
\nu\lambda_{d}t^{2-d}-\lambda_d\nu^2\sum\limits_{k=1}^n
\nu_k^2r(B,\mathbf x^*_k)^{2-d}+ \nu^2 E(X^*,\Delta,B) +o(1)
\leq\\
\nu\lambda_{d}t^{2-d}-\lambda_d\nu^2\sum\limits_{k=1}^n
\nu_k^2r(B,\mathbf x_k)^{2-d}+ \nu^2 E(X,\Delta,B) +o(1),\ t\to 0,
\end{multline}
where $\nu_k=\delta_k$,
$\nu=\left(\sum\limits_{k=1}^n|\delta_k|\right)^{-1}$,
$k=1,\ldots,n$. As $B$ is the rotation domain, harmonic radii
$r(B,\mathbf x)$ take equal values at all points  $\mathbf x$,
lying on one circle from $\Omega$.  Hence,
$$
\sum\limits_{k=1}^n \nu_k^2r(B,\mathbf
x^*_k)^{2-d}=\sum\limits_{k=1}^n \nu_k^2r(B,\mathbf x_k)^{2-d}
$$ and it remains to take the limit as $t\to 0$ in \eqref{as2} to complete the proof of Theorem~1.

Let the domain $B$ and the collections $X=\{\mathbf
x_k\}_{k=1}^n$, $X^*=\{\mathbf x_k^*\}_{k=1}^n$,
$\Delta=\{\delta_k\}_{k=1}^n$ be as in Theorem~2.  We can assume
that the boundary $\partial B$ represents a continuously
differentiable surface in $\rr$.  Put
$\sigma_k=\text{sgn}\,\delta_k$, $\mu_k=|\delta_k|^{1/(d-2)}$,
$k=1,\ldots,n$, $\Lambda =\{\sigma_k\}_{k=1}^n$,
$\Psi=\{\mu_k\}_{k=1}^n$. The condenser $C(t;B,X,\Delta,\Psi)$
admits \emph{a potential function} $u$, which is continuous in
$\overline B$, harmonic in \linebreak $B\setminus\left(E(\mathbf
x_1,\mu_1{t})\cup\ldots\cup E(\mathbf{x_n}, \mu_n{t})\right)$ and
$u=0$ on $\partial{B}$, $u=\sigma_k$ on the set $E(\mathbf x_k,
\mu_k t)$, $k=1,\ldots,n$, \cite{Lan}. Denote by $I$ the set of
points in $B$, where $u=0$, and by $D_l$, $l=1,\ldots,q$, the
connected components  $B\setminus I$. Suppose $X_l$ is the set of
all points from $X$, lying in the domain $D_l$.

For the collection of points $X_l=\{\mathbf y_{sl}\}_{s=1}^{N_l}$
let us define $\Lambda_l=\{\sigma_{sl}\}_{s=1}^{N_l}$ and
$\Psi_l=\{\mu_{sl}\}_{s=1}^{N_l}$ according to the rule
$\sigma_{sl}=\sigma_p$, $\mu_{sl}=\mu_p$, if $\mathbf
y_{sl}=\mathbf x_p$, $\mathbf x_p\in X$. According to the
Dirichlet principle \cite{Sob}
\begin{multline}\label{cap} {\mathrm
{cap}}\,C(t;B,X,\Lambda,\Psi)=\int\limits_B |\nabla u|^2\:dx
=\sum\limits_{l=1}^q \int\limits_{D_l} |\nabla u|^2\:dx
=\sum\limits_{l=1}^q{\mathrm
{cap}}\,C(t;D_l,X_l,\Lambda_l,\Psi_l).
\end{multline}

Note that the points lying in the domain $D_l$ have the same
charge (either $1$ or $-1$) so that according to
Lemma~\ref{lemma2} the following equality holds:
\begin{equation}\label{cap1}
{\mathrm {cap}}\,C(t;B,X,\Lambda,\Psi)=\sum\limits_{l=1}^q{M_l},
\end{equation}
where $M_l=M(\Gamma(t;D_l,X_l,\Psi_l))$. Let $n_l$ denote the
number of half-planes ${L}_j$, containing at least one point from
$X_l$, $0\leq j\leq m-1$. It is clear that
$\sum\limits_{l=1}^{q}n_l\geq m$.  Convexity of the function $1/x$
implies that for any positive numbers $v_l$, $\alpha_l$,
$\sum\limits_{l=1}^q\alpha_l=1$, $l=1,\ldots,q$, the following
inequality holds:
$$
\left(\sum\limits_{l=1}^q \alpha_l(\alpha_l^{-1}
v_l)\right)^{-1}\leq\sum\limits_{l=1}^q \alpha_l(\alpha_l^{-1}
v_l)^{-1},
$$
or
\begin{equation}\label{in1}
\left(\sum\limits_{l=1}^q v_l\right)^{-1}\leq \sum\limits_{l=1}^q
\alpha_l^2 v_l^{-1}.
\end{equation}
Obviously, inequality \eqref{in1} remains valid for any
non-negative  $\alpha_l$,  $\sum\limits_{l=1}^q\alpha_l\geq 1$
and, moreover,
\begin{equation}\label{in2}
\left(\sum\limits_{l=1}^q v_l\right)^{-1}\leq
\frac{1}{q^2}\sum\limits_{l=1}^q v_l^{-1}.
\end{equation}
Then it follows from  (\ref{cap1}) and \eqref{in1} that
\begin{equation}\label{mod1}
|C(t;B,X,\Delta,\Psi)|=\left(\sum\limits_{l=1}^q
M_l\right)^{-1}\leq \sum\limits_{l=1}^q\frac{n_l^2}{m^2} M_l^{-1}.
\end{equation}

Denote by $\Gamma_{lj}^+$ the family of curves from
$\Gamma(t;D_l,X_l,\Psi_l)$ lying in the dihedral angle
$\{(r,\theta,\mathbf x'): \theta_j\leq\theta\leq\theta_{j+1}\}$
and by $\Gamma_{lj}^-$ the family of curves from
$\Gamma(t;D_l,X_l,\Psi_l)$ lying in the dihedral angle
$\{(r,\theta,\mathbf x'): \theta_{j-1}\leq\theta\leq\theta_{j}\}$,
$\theta_m=\theta_0+2\pi,$ $\theta_{-1}=\theta_{m-1}-2\pi$,
$l=1,\ldots,q,$ $j=0,\ldots m-1$. According to the property
\eqref{prop1} we have
\begin{equation}\label{in3}
M_l\geq \sum\limits_{j=0}^{m-1}\hspace{-0.05cm}\left.
^{'}\hspace{-0.1cm}\left(M(\Gamma_{lj}^-)+M(\Gamma_{lj}^+)\right)\right.,
\end{equation}
where the prime at the summation sign means that the summation is
taken over those indices $j$, $j=0,\ldots,m-1,$ for which the
half-plane $L_j$ contains at least one point from $X_l$. Note that
the total number of terms in this sum equals $2n_l$. Inequalities
\eqref{in3} and \eqref{in1} imply that
\begin{multline}\label{mod2}
\sum\limits_{l=1}^q\frac{n_l^2}{m^2} M_l^{-1}\leq
\sum\limits_{l=1}^q\frac{4n_l^2}{4m^2} \left(
\sum\limits_{j=0}^{m-1}\hspace{-0.05cm}\left.
^{'}\hspace{-0.1cm}\left(M(\Gamma_{lj}^-)+M(\Gamma_{lj}^+)\right)\right)^{-1}\right.\leq
\\ \sum\limits_{l=1}^q\frac{1}{4m^2}
\sum\limits_{j=0}^{m-1}\hspace{-0.05cm}\left.
^{'}\hspace{-0.1cm}\left(M(\Gamma_{lj}^-)^{-1}+M(\Gamma_{lj}^+)^{-1}\right).\right.
\end{multline}

Next, consider the symmetric configuration that is the condenser
$C(t;B,X^*,\Lambda,\Psi)$. Let $X^*_0$ be the collection of points
in $X^*$, lying on the half-plane $\{(r,\theta,x'):\theta=0\},$
$X_0^*=\{\mathbf y_{s}^*\}_{s=1}^{K}$.   If $\mathbf y_s^*=\mathbf
x_p^*,$ $\mathbf x_p^*\in X^*$, then we define ${\mu}_{s}^*$ by
the equality $\mu_{s}^*=\mu_p$.  In view of the symmetry of the
condenser $C(t;B,X^*,\Lambda,\Psi)$ we have
$$
{\rm cap}\,C(t;B,X^*,\Lambda,\Psi)=m{\rm cap}\,C(t;B\cap
P_0,X_0^*,\Lambda_0^*,\Psi_0^*),
$$
where $P_0=\{(r,\theta,\mathbf x'): -\pi/m<\theta<\pi/m\}$,
$\Lambda_0^*=\{\sigma_s^*\}_{s=1}^{K}$,
$\sigma_1^*=\ldots=\sigma_K^*=1$,
$\Psi_0^*=\{\mu_{s}^*\}_{s=1}^{K}$. Using the symmetry of the
condenser $C(t;B\cap P_0,X_0^*,\Lambda_0^*,\Psi_0^*)$ and
Lemma~\ref{lemma2}, we conclude that
$$
{\rm cap}\,C(t;B,X^*,\Lambda,\Psi)=2m M(\Gamma_0^*),
$$
where $\Gamma_0^*$ is the family of those curves from the family
$\Gamma(t;B\cap P_0,X_0^*,\Psi_0^*)$ that lie in the angle
$\{(r,\theta,x'): 0\leq\theta\leq\pi/m\}$. Write
$\phi_k(\mathbf{x})$ for the reflection with respect to the
hyper-plane $\{(r,\theta,\mathbf x'):\theta=\pi k/m\ \text{or}\
\theta=\pi+\pi k/m\}, \ k=1,\ldots,2m-1$. For each curve
$\gamma_0^*\in \Gamma_0^*$ define
$\gamma_k^*=\phi_k(\gamma_{k-1}^*)$ and the curve
$\gamma^*=\cup_{k=0}^{2m-1}\gamma_k^*$. In other words, $\gamma^*$
is the curve symmetric with respect to the group $\Phi$ (see
section~2) and comprising $2m$ consecutive reflections
$\gamma_0^*$. Let $\Gamma^*$ be the family of curves $\gamma^*$.
The composition principle and the symmetry of the family
$\Gamma^*$ \cite[c.21]{Alf}, \cite[c.178,179]{Fuglede}  imply that
$$
M(\Gamma_0^*)=2m M(\Gamma^*).
$$
Hence,
\begin{equation}\label{eq:cap2}
|C(t;B,X^*,\Lambda,\Psi)|=\frac{M(\Gamma^*)^{-1}}{4m^2}.
\end{equation}

We now apply dissymmetrization described in Lemma~\ref{LemmaDis}.
As dissymmetrization induces a metric in each direction which
preserves length and volume, we have
$$
M(\Gamma^*)=M({\rm Dis}\,\Gamma^*),
$$
where ${\rm Dis}\,\Gamma^*=\{{\rm Dis}\,\gamma^*:\gamma^*\in
\Gamma^*\}$. It is easy to see from the construction of
dissymmetrization (details can be found in
\cite[pp.63-64]{KalmPrilEn}) that the family
$\text{Dis}\,\Gamma^*$ is minorized by the separated
$\Gamma_{lj}^-$ and  $\Gamma_{lj}^+$. Due to the property
\eqref{prop2} and Lemma~\ref{lemma4}, we then have
$$
\sum\limits_{l=1}^q\frac{1}{4m^2}\sum\limits_{j=0}^{m-1}\hspace{-0.05cm}\left.
^{'}\hspace{-0.1cm}\left(M(\Gamma_{lj}^-)^{-1}+M(\Gamma_{lj}^+)^{-1}\right)\right.\leq
\frac{1}{4m^2} M(\text{Dis}\, \Gamma^*)^{-1}=\frac{1}{4m^2}
M(\Gamma^*)^{-1}.
$$
In view of  (\ref{mod1}), (\ref{mod2}), (\ref{eq:cap2}) the above
inequality leads to
$$
|C(t;B,X,\Lambda,\Psi)|\leq |C(t;B,X^*,\Lambda,\Psi)|.
$$
It remains to apply the asymptotic formula \eqref{eq:asimptotika}
for the condenser modulus following the same line of argument as
in the proof of Theorem~1. This completes the proof of Theorem~2.

\textbf{Funding:} This work was supported by the Russian Basic
Research Fund   [grant number 20-01-00018].

\end{document}